\documentclass[graybox]{svmult}

\usepackage{times}
\usepackage{amsmath}
\usepackage{amssymb}
\usepackage{amscd}
\usepackage{stmaryrd}

\usepackage{xspace}
\input xy
\xyoption{all}

\usepackage{helvet}         % selects Helvetica as sans-serif font
\usepackage{courier}        % selects Courier as typewriter font
\usepackage{type1cm}        % activate if the above 3 fonts are
\usepackage{makeidx}         % allows index generation
\usepackage{graphicx}        % standard LaTeX graphics tool
\usepackage{multicol}        % used for the two-column index
\usepackage[bottom]{footmisc}% places footnotes at page bottom

\input xy
\xyoption{all}
\xyoption{knot}

\def\S{{\mathfrak S}}

\def \la {\left\langle}
\def \ra {\right\rangle}

\def\hr{{{\mathcal H}_r^q}}
\def\hh{{{\mathcal H}_2^q}}
\def\hhh{{{\mathcal H}_3^q}}
\def\h{{{\mathcal H}^q}}
\def\hdia{{{\h}^\Delta}}

\def\PP{{\mathfrak{PP}}}

\def\pK{{\mathfrak{L}_q^\Delta}}

\def\dia{{\C[GL_q(V)]^\Delta}}
\def\T{{\mathbb T}}
\def\Uq{{U_q{\mathfrak{gl}}\,(V)}}
\def\C{{\mathbb C}}
\def\N{{\mathbb N}}

\def\S{{\mathfrak S}}
\def\I{{\mathfrak I}}

\def\calh{{\cal H}}

\newcommand{\beq}{\begin{equation}}
\newcommand{\eeq}{\end{equation}}
\newcommand{\ba}{\begin{array}}
\newcommand{\ea}{\end{array}}
\newcommand{\beqa}{\begin{eqnarray}}
\newcommand{\eeqa}{\end{eqnarray}}

\newcommand{\g}{\mathfrak {g}}

\newcommand{\Vt}[1]{V^{\otimes #1}}
\newcommand{\Vst}[1]{V^{\ast \otimes #1}}

\begin{document}

\title*{Quantum Diagonal Algebra and \\Pseudo-Plactic Algebra }
\author{Todor Popov}
\institute{  INRNE, Bulgarian Academy of Sciences, 72 Tsarigradsko Chaussee, 1784 Sofia, Bulgaria\\
\& American University in Bulgaria
%\institute{Todor Popov \at Name, Address of Institute, 
\email{tpopov@aubg.edu}
%\and Name of Second Author \at Name, Address of Institute \email{name@email.address}
}

\maketitle

\abstract{The subalgebra of diagonal elements of a
quantum matrix group has been conjectured 
by Daniel Krob and Jean-Yves Thibon to be isomorphic
to a cubic algebra, coined the quantum pseudo-plactic algebra.
We present a functorial approach to the conjecture
through the quantum Schur-Weyl duality between the quantum group and the Hecke algebra. 
The relations of the  quantum  diagonal subalgebra 
are found to be the image of the  braid relations of the
underlying  Hecke algebra by an appropriate 
Schur functor
%through an appropriate Schur functor
which gives a straightforward proof of the conjecture.
 }

\section{Introduction}
\label{sec:1}

The diagonal  elements %of the coordinate ring
of a quantum matrix  group $\C[GL_q(V)]$ close a subalgebra which is the noncommutative avatar of the  algebra of the functions on the torus. The resulting
 quantum diagonal algebra $\C[GL_q(V)]^\Delta \subset \C[GL_q(V)]$ 
provides
a noncommutative character theory of quantum group comodules
which is a lifting of the commutative symmetric functions.
%Let us denote by $\T$ the torus and 
We identify the functions on the quantum torus $\T$ with the
 subspace of $End(V)^\ast$ stable  by  
 the transposition $\tau$,  $\tau(x^i_j)= x^j_i$.
Krob and Thibon conjectured \cite{KT}  that the algebra $\C[GL_q(V)]^\Delta$
spanned by $x^i_i\in \C[GL_q(V)]$ is isomorphic to  the {\it quantum pseudo-plactic algebra} defined as the quotient 
$\PP_q(\T) \cong \C (q)\la \T \ra / (\pK(\T))
				$
of the free diagonal algebra $\C(q)\la \T \ra$ by the  
ideal $({\pK}(\T))$ generated by
\beq
\label{ppK}
\ba{rclcrcl}
         {\pK}^{i_1,i_2}_{i_3}  &:=&
   [  [ x_{i_1}^{i_1},x_{i_3}^{i_3} ], x_{i_2}^{i_2}
]      &\qquad&\mbox{with}&\qquad
&i_1<i_2<i_3
   \\ [4pt]
 {\pK}^{i_1,i_1}_{i_2}  &:=&
        [[ x_{i_1}^{i_1}, x_{i_2}^{i_2} ], x_{i_1}^{i_1} ]_{q^2}
& & \mbox{with}& & i_1<i_2
 \\[4pt]
 {\pK}^{i_1,i_2}_{i_2}  &:=&
  [  x_{i_2}^{i_2}, [ x_{i_1}^{i_1}, x_{i_2}^{i_2} 
]]_{{q^2}}               && \mbox{with} & & i_1<i_2
         \ea \ . \eeq 
				%\[
				%\PP_q(\T) \cong \C (q)\la \T \ra / (\pK(\T))
				%\]
		Stated differently, the ideal of relations of the diagonal subalgebra $\C[GL_q(V)]^\Delta$ is generated by the cubic relations (\ref{ppK}) and there are no 
		relations in higher order
		which are independent from the cubic ones, eq (\ref{ppK}).

				We introduce  a Schur bifunctor %in chapter \ref{sch} which maps
				the tower of Hecke algebras $\hr$ into 
				the coordinate ring of
				the quantum group $\C[GL_q(V)]$. It maps $(\hr,\hr)$-modules into $(\Uq,\Uq)$-modules. The polarization functor denoted by $\xymatrix{\ar@{~>}[r]&}$ is the  adjoint functor of the Schur bifunctor.
				The polarization of the diagonal subalgebra 
				$\C[GL_q(V)]^\Delta$ is defined to be the {\it diagonal Hecke algebra ${\h}^\Delta$}	
%$\PP_q(\T)$
%\beq
%\label{conj}
%\C[GL_q(V)]^\Delta \cong \PP_q(\T) \ .
%\eeq
\beq
\label{conj}
\xymatrix{ 
\PP_q(\T) \ar@{~>}[d]&\cong &
\C[GL_q(V)]^\Delta \ar@{~>}[d] \ar@{^{(}->}[r]  &  
 \C[GL_q(V)] =\bigoplus_{r\geq 0} \C[GL_q(V)]_r  \ar@{~>}[d] \\ 
 \PP_q  &  \cong & \h^\Delta %=\bigoplus_{r\geq 0} \hr 
\,\, \, \,\,\, \vspace{8pt} \ar@{^{(}->}[r]   &
\vspace{8pt}\h =\bigoplus_{r\geq 0} \hr 
\ .}
\eeq
The quantum Weyl action \cite{LS} is stabilizing 
the diagonal $\C[GL_q(V)]_r^\Delta$
in $\C[GL_q(V)]_r$.
It has a counterpart,  an ``adjoint" $\hr$-action 
on $\hr^\Delta$. We consider also the polarization
$\PP_q$ of the 
quantum pseudo-plactic algebra  $\PP_q(\T)$ which we refer to as
the pre-plactic algebra defined as a factor algebra.
We prove that the ideal of $\PP_q$ induced by 
the unique polarized pseudo-plactic relation, namely
\beq
\label{prepl}
\pK= [  [ x_{1}^{1},x_{3}^{3} ], x_{2}^{2}] \quad \Leftrightarrow \quad
\pK=(x_1^1 x_3^3 x_2^2 - x_3^3 x_1^1 x_2^2) -
(x_2^2 x_1^1 x_3^3 - x_2^2 x_3^3 x_1^1)
\eeq
%$$[  [ x_{1}^{1},x_{3}^{3} ], x_{2}^{2}]=0$$ 
contains all relations of the diagonal Hecke algebra 
${\h}^\Delta$, that is,
one has the isomorphism  of $\h$-modules between
the pre-plactic algebra and the diagonal Hecke algebra
$$\PP_q \cong \h^\Delta \ .$$

The pre-plactic relation (\ref{prepl})
 is the difference of the Knuth relations
of the plactic monoid \cite{Plaxique} therefore
the pre-plactic algebra   $\PP_q$ is a lifting of the Poirier-Reutenauer  algebra \cite{PR}. 
%From our perspective
%the polarization of the Lascoux and Schtuzenberger
%its relation is the difference of the Knuth relations
%of the plactic monoid
%\[
%x_1^1 x_3^3 x_2^2 - x_3^3 x_1^1 x_2^2 = 
%x_2^2 x_1^1 x_3^3 - x_2^2 x_3^3 x_1^1
%\]
 %It follows that the Poirier-Reutenauer algebra (\cite{PR})
%is  a quotient of the pre-plactic algebra $\PP_q$.
Both Poirier-Reutenauer  $\mathfrak{PR}$  and
 pre-plactic algebra $\PP_q$ \cite{TP2}
%have   Hopf algebra structure 
are Hopf algebra
 quotients of the Malvenuto-Reutenauer algebra
 $\mathfrak{MR}$ \cite{MR}
\[ 
\mathfrak{PR} \subset \PP_q \subset  \mathfrak{MR} \ .
\]
%There exists a Schur functor 
% mapping $\mathfrak{PR}$ 
%to the algebra of the plactic monoid,
There exists also a Schur functor 
 %$q$-deformation
 $PS_q$ \cite{LP3} mapping the 
$q$-deformation of $\mathfrak{PR}$ 
to  $q$-deformation $PS_q(V)$ of the plactic algebra
(deformed parastatistics algebra \cite{D-VP,LP2}).

In  Chapter \ref{pre-pl} by applying a Schur bifunctor to
the $\h$-modules
$\PP_q \cong \h^\Delta$
we prove the conjecture of Krob and Thibon
\[
\PP_q(\T) \cong \C[GL_q(V)]^\Delta \ .
\]
The cubic relations ${\pK(\T)}=0$ are playing a role similar to the Knuth relations of the plactic algebra in the theory of noncommutative symmetric functions. Moreover
the
pre-plactic relation $\pK=0$ acquire in the process of our proof
 clear geometrical meaning,
it is nothing but the 
braid relation of the Hecke algebra. 
%\[
%\T^\ast=\{ x\in {\rm End(V)}^\ast | x=x^\tau \} 
%\]
 %(or the restriction of
%$V \otimes V^\ast$ to the diagonal 
%see  below).

\section{Schur functor and its adjoint}
\label{sch}
%{Quantum Schur-Weyl duality}

The algebra of functions on the quantum group 
$\C[GL_q(V)]$ coacts on itself, it is a naturally
a $(\C[GL_q(V)],\C[GL_q(V)])$-comodule. The ``regular representation" of $\C[GL_q(V)]$ according to 
the Peter-Weyl theorem
has a decomposition into  a product of 
 left and right irreducible $[GL_q(V)]$-comodules
$$\C[GL_q(V)]_r\cong 
 %(V^\ast)^{\otimes r} \otimes_{\hr} V^{\otimes r}
%\cong
\bigoplus_{\lambda \vdash r} S_{\lambda}(V^{\ast}) 
\otimes%_{\hr} 
S^{\lambda} (V)\ .$$
In their seminal paper Faddeev-Reshetikhin-Takhtajan \cite{FRT}   defined the algebra of functions on 
the general linear quantum group $\C[GL_q(V)]$
as the commutant of the action
of the Hecke algebra ${{\mathcal H}}_2(q)$.
The Hecke algebra ${{\mathcal H}}_2(q)$ is
represented by  the Drinfeld-Jimbo quantum $R$-matrix $\hat{R}_q\in End(\Vt 2)$
\beq
\label{frt}
\C[GL_q(V)] = \C(q)\la W^\ast \ra /(\hat{R}_q W\otimes W- W \otimes W \hat{R}_q) \quad W =End(V)^\ast\cong V^\ast \otimes V \eeq
In other words the homogeneous elements of $\C[GL_q(V)]$
are the coinvariants of the natural $\hr$-action
\[
\C[GL_q(V)]_r\cong (W^{\otimes r})^{\hr}\cong 
 (V^\ast)^{\otimes r} \otimes_{\hr} V^{\otimes r} \ .
\]
Here the decorated tensor product  $\otimes_{\hr}$
stands for the quotient relating the right and left 
$\hr$-action,  e.g. $ (W^{\otimes 2})^{{\mathcal H}^q_2} :=
W^{\otimes 2}/( W \otimes W \hat{R}_q -
\hat{R}_q W\otimes W)$.

By duality the two-sided comodules $\C[GL_q(V)]_r$ are 
bimodules of the
quantum universal enveloping algebra 
$U_q(\g)$. The quantum Schur-Weyl duality is the double commutant 
property
of the action of the quantum universal enveloping algebra 
$\Uq$
and the action of the  Hecke algebra $\hr$. The Schur-Weyl duality allows to build the Schur functor 
which maps the category of representations
of the Hecke algebra $\hr-mod$
to the category of representations of the 
quantum universal enveloping algebra $\Uq-mod$.
%\[
%{\rm End}_{\Uq} \Vt r = \hr \qquad \qquad
%{\rm End}_{\hr} \Vt r = \Uq
%\]

%The irreducubles in  $\hr-mod$ are constructed as
{\bf Orthogonal idempotents in $\hr$.} An orthogonal
idempotent $e_{\lambda}(T)$ in $\hr$ is parametrized by a partition
$\lambda$ of $r$, $\lambda \vdash r$ and a Standard Young Tableau $T$
with shape $\lambda$, $T\in STab(\lambda)$. 
Different idempotents are orthogonal
$$
e_\lambda (T) e_\mu(T')=e_\lambda (T) \delta_{\lambda \mu}\delta_{TT'} \ .
$$
A system of orthogonal idempotents provides a partition of unity
$$
\sum_{\lambda \vdash r} \sum_{sh(T)=\lambda} e_\lambda (T)
= 1 \!\! 1_{\hr} \ .
$$

An  irreducible right $\hr$-module 
is constructed as the ideal given by multiplication with 
an idempotent $e_\lambda$ from the left, 
$S^\lambda =  e_\lambda \hr $.  Idempotents $e_\lambda(T)$ and
$e_\lambda(T^\prime)$
with different Young Tableaux $T,T^\prime \in STab(\lambda)$
having same shape  $\lambda$ lead to isomorphic
$\hr$-modules $S^{\lambda} \cong S^{\lambda(T)} 
\cong S^{\lambda(T^\prime)}$,
so we often suppress $T$.

The Schur functor 
$S^\lambda (V)= S^\lambda \otimes_{\hr} \Vt r$
maps the irreducible right $\hr$-module
$S^\lambda = e_\lambda \hr$ into an irreducible 
right $\Uq$-module. Similarly on defines
and  irreducible left $\hr$-module
$S_\lambda =  \hr e_\lambda$ and
irreducible left $\Uq$-module
$S_\lambda (V^\ast)= \Vst r \otimes_{\hr} S_\lambda  $.
The polarization is a functor adjoint to the Schur functor,
it maps a $\Uq$-module $S^\lambda(V)$ with $|\lambda|=r$ into its underlying
$\hr$-module $S^\lambda$.

%Beside the grading by the degree $r$
The coordinate ring is $\N^d$-bigraded % $d=\dim V$ 
by the weight{\footnote{The weight
of the multi-index $A$ is the $d$-dimensional vector $(w_1(A), \ldots , w_d(A))$
 %$(w_i(A))_{1\leq i \leq d}$
defined as $w_i(A)= \#\{ a_k=i|1 \leq  k\leq r\}$} of the multi-indices $A$ and $B$
}
$$ \C%_{\hat{R}_q}
[GL_q(V)]_{r\, B}^{\, A }\cong
%(V^{\ast \otimes r})^A\otimes_{\hr} (V^{\otimes r})_B
%\cong
\bigoplus_{\lambda \vdash r} S^A_{\lambda}(V^{\ast}) \otimes_{\hr} S_B^{\lambda} (V) \qquad A,B\in \{1, \ldots , d=\dim V \}^r \ .$$
The polarization of $[GL_q(V)]_r$ is 
its component %with multi-indices $\alpha$ and $\beta$ of 
of weight $1^r$ (we suppose that we have chosen $V$ such that $\dim V=r$).
Equivalently the multi-indices  $\alpha$ and $\beta$
of weight $1^r$  are words
of permutations $\alpha, \beta\in \S_r$.
The polarization of the $(\Uq,\Uq)$-module $\C[GL_q(V)]_{r}$ is isomorphic to the $(\hr,\hr)$-bimodule $\hr$
yielding the decomposition % into irreducibles Specht module
 of the regular representation of $\hr$ 
  %underlying   $\C[GL_q(V)]_{r}$ 
\beq
\label{coset}
(\hr)^\alpha_\beta:=\C[GL_q(V)]_{r\, \beta}^{\, \alpha } =
 \bigoplus_{\lambda \vdash r} 
 S_\lambda^\alpha \otimes_{\hr} S^\lambda_\beta \  \ .
\eeq
Every left (right) $\hr$-module $S_\lambda$ ($S^\lambda$) appears in the regular representation with multiplicity equal to its dimension
$f_\lambda= \dim S_\lambda$.

 For generic $q$
the regular representation $\hr$ is 
 isomorphic to its specialization at $q=1$, {\it i.e.},
to the regular representation of $\C[\S_r]$. Hence
$\C[\S_r]$ can be seen as the polarization of
the commutative algebra $\C[GL(V)]$.
A permutation $\alpha \in {\S}_r$ is given
by the two row bijective correspondence (with commuting biletters) or equivalenty
by its word
\[
\alpha = 
\left(\ba{ccc} 
1 & \ldots & r 
\\
\alpha_1 & \ldots & \alpha_r
\ea \right) \qquad \alpha = (\alpha_1 \ldots \alpha_r) \ .
\]
The inverse permutation $\alpha^{-1}$ is simply obtained
by exchanging the two rows
\[
\alpha^{-1} = 
\left(\ba{ccc} 
\alpha_1 & \ldots & \alpha_r \\
1 & \ldots & r 
\ea \right)= 
\left(\ba{ccc} 
1 & \ldots & r\\
\alpha_1^{-1} & \ldots & \alpha_r^{-1}
\ea \right) %\qquad \alpha = (\alpha_1 \ldots \alpha_r) \ .
\]
in which we have rearranged the commuting biletters.
 %$\left( \ba{c} i \\ \alpha_i \ea \right)$.
More generally any two row bijection yields a permutation
representable as a product 
\beq
\label{coset1}
\left(
\ba{c} \alpha \\ \beta
\ea \right) = \left(
\ba{c} \alpha_1 \ldots \alpha_r \\ 1 \ldots r
\ea \right) \left(
\ba{c}  1 \ldots r
 \\ \beta_1 \ldots \beta_r \ea \right) = \alpha^{-1} \beta \ .
\eeq
The  permutation $\alpha^{-1} \beta$ can be thought as an element in the double coset
where the right $\C[\S_r]$-action is by place permutation
and the left $\C[\S_r]$-action is by substitution
$$ \alpha^{-1} \beta= \left(
\ba{c} \alpha \\ \beta
\ea \right)  \in \C[\S_r] \cong \C[\S_r] \otimes_{\C[\S_r]} \C[\S_r] \ . $$
The commutativity of the biletters is expressed
by the coset notation $\otimes_{\C[\S_r]}$.

We now come back to the Hecke algebra $\hr$
isomorphic to the permutation group algebra $\hr \cong \C[\S_r]$ for generic $q$.  It has a basis $T_{\sigma}\in \hr$ indexed by the permutations $\sigma \in S_r$. 
We intorduce another basis $T^{\sigma}:=T_{\sigma^{-1}}$

In parallel with $\C[\S_r]$
one has the double coset for $(\hr)^\alpha_\beta$,
cf. eq. (\ref{coset}) with basis 
\beq
\label{cs}
T^\alpha_\beta \in \hr \cong \hr \otimes_{\hr}\hr \qquad  
T^\alpha_\beta:= T^\alpha \otimes_{\hr} T_\beta
 %=1\!\!1 \otimes_{\hh}
%T_{\alpha^{-1}} T_\beta  
=T_{\alpha^{-1}} \otimes_{\hh} T_\beta  
\qquad \alpha, \beta \in \S_r \ .
\eeq
%in parallel with eq. (\ref{coset1}).

{\bf Polarization of $\C[GL_q(V)]$ relations. }
% We now consider as an example the polarization of 
%the ideal
%$\hat{R}_q W\otimes W - W \otimes W\hat{R}_q$.
By evaluation  of the Drinfeld-Jimbo $R$-matrix $\hat{R}_q$
with indices $i,j$ running in the range $ 1 \leq i,j \leq d=\dim V$
\beq
\label{DJ}
\hat{R}_q=% q \sum_{i} e_{i}^{i} \otimes  e_{i}^{i} +
\sum_{i,j}  q^{\delta_{ij}}e_{j}^{i} \otimes  e_{i}^{j} +
 (q-q^{-1} )\sum_{ i<j} e_{j}^{i} \otimes  e_{j}^{i}
\quad \qquad e_{j}^{i} \in \mathfrak{gl}(V)  \ .
\eeq
%eq. (\ref{DJ}) % in  the FRT relations
we get the Faddeev-Reshetikhin-Takhtajan relations \cite{FRT}
$\hat{R}_q W\otimes W = W \otimes W\hat{R}_q$ 
of the coordinate ring of the quantum group $\C[GL_q(V)]$
\beq
\label{frt2}
 \ba{lclcl}
 x^j_k x^i_k =q  x^i_k x^j_k  &\quad&
 x_j^k x_i^k =q  x_i^k x_j^k & \quad \qquad&j > i  \\
  x^j_l x^i_k =x^i_k x^j_l +
 (q- q^{-1})x^i_l x^j_k && x^j_k x^i_l = x^i_l x^j_k & & j > i \quad l > k
 \ea \ .
%\label{qg}
 \eeq
These relations span an ideal which is also a coideal for the  coaction 
$\Delta x^i_j = \sum_{k} x^i_{k} \otimes x^k_j$.
By duality it is a $(\Uq,\Uq)$-module.
The polarization of the ideal generators $\hat{R}_q W\otimes W - W \otimes W\hat{R}_q$, that is, the weight $1^2$ relations
of $\C[GL_q(2)]$,
 yields a
%Its polarization it the Hecke algebra 
$({\mathcal H}^q_2,{\mathcal H}^q_2)$-module ${\mathcal H}^q_2$
%spanned by
\beq
\label{pol}
T^{21}_{\,21} =T^{12}_{\,12} + (q-q^{-1})T^{12}_{\,21} \qquad
T^{12}_{\,21}=T^{21}_{\,12} \ .
\eeq
With the help of the identification (\ref{cs})
one has $T^{12}_{\,12} = % T_{12}=
1\!\! 1  \otimes 1\!\! 1 $
\[ 
T^{21}_{\,21}= T_{s_1} \otimes_{{\mathcal H}^q_2}
T_{s_1} = 1\!\! 1 \otimes_\hh (T_{s_1})^2
\qquad 
T^{12}_{\,21} = 1\!\! 1 \otimes_\hh T_{s_1} %= T_{21}
=T_{s_1}\otimes_\hh 1\!\! 1 = T^{21}_{12} \ .
\]
Hence the polarization (\ref{pol})
of eq. (\ref{frt2})
is equivalent to  the two factorizations of the Hecke relation
(taken together)
%The basis of the two-dimensional Hecke bimodule ${\mathcal H}^q_2$ is $T_{12}=1\!\! 1$, $T_{21}=T_{s_1}$.
\beq
\label{hr}
(q^{-1} 1\!\! 1 + T_{s_1}) \otimes_\hh (q 1\!\! 1 - T_{s_1})
=0=(q 1\!\! 1 - T_{s_1})\otimes_\hh(q^{-1} 1\!\! 1 + T_{s_1}) \ 
\eeq
Hence by applying the polarization functor 
to the quantum matrix group $\C[GL_q(V)]$
relations
$\hat{R}_q W\otimes W - W \otimes W\hat{R}_q$
 %the quadratic relations
one gets the relations of the Hecke algebra $\hh$.

Conversely, the normalization of %the ideal of $\hh$ generated by
the factorization of the Hecke relations, cf. eq. (\ref{hr}) (after division
by $[2]=q+ q^{-1}\neq 0$ for $q^2\neq -1$)
leads to  
%the product of the 
two orthogonal idempotents in $\hh$ 
%( one has the 
%unit partition $ 1\!\! 1_{\hh}=e_{2}+e_{1^2}$)
\[
e_{2}e_{1^2}= 0 = e_{1^2} e_{2} 
\qquad \qquad (1\!\! 1_{\hh}=e_{2}+e_{1^2})
\]
where the $q$-symmetrizer  and the $q$-antisymmetrizer are respectively
$$
 e_{2}:=%E^{(2)}(q)=
 \frac{1}{[2]}
\left(  q^{-1} 1\!\! 1 + T_{s_{1}}\right) \ ,\qquad
e_{1^2}:=%E^{(1^2)}(q)=
 \frac{1}{[2]}\left( q 1\!\! 1 - T_{s_{1}} \right) \ .
$$
%The  Schur functor maps the right $\hh$-module 
%$S_{\lambda}=  e_{\lambda} \hh$ to 
 %the  right $\Uq$-module $S_{\lambda}(V)$, e.g.,
%$
%S^{2}(V)= S^{2} \otimes_{\hh} V^{\otimes 2}=
%e_{2} V^{\otimes 2}$  and $S^{1^2}(V)= e_{1^2} V^{\otimes 2}  $. Similarly of the left 
%$\Uq$-modules $S_{\lambda}(V^\ast)$.

%The coordinate ring 
Half of the relations 
of the quantum matrix group $\C[GL_q(V)]$, cf. eq. (\ref{frt})
are obtained through the Schur bifunctors
$$e_{2} W\otimes W e_{1^2}=S_{2}(V^\ast) \otimes S^{1^2}(V)=
V^{ \ast\otimes  2} \otimes_{\hh}  S_{2} \otimes  S^{1^2} \otimes_{\hh} \Vt 2 $$
where $S^{1^2}(V)= e_{1^2} V^{\otimes 2} =
S^{1^2} \otimes_{\hh} V^{\otimes 2} $ and
$S_{2}(V^\ast) = V^{ \ast \otimes 2} e_{2} = 
V^{\otimes 2} \otimes_{\hh} S^{2}$.
The $\hh$-action $\rho$ of the projectors $e_{2}$ and $e_{1^2}$
is throught the multiplication by the R-matrix eq.(\ref{DJ}),
$\rho(T_{s_1})= \hat{R}_q \in {\rm End}(\Vt 2)$.
The other half of the relations in eq. (\ref{frt})  
, (the `` missing relation'' after Yuri Manin \cite{Manin})
$$ 
e_{1^2} W\otimes W e_{2}=S_{1^2}(V^\ast) \otimes S^{2}(V)=
V^{ \ast\otimes  2} \otimes_{\hh}  S_{1^2} \otimes  S^{2} \otimes_{\hh} \Vt 2 \ .$$

The above relations define the so called
left and right quantum semi-groups \cite{Manin}, respectively.
Taken together they 
span the ideal of the quantum group relations eq. (\ref{frt}). 
Equivalently one has
the  short exact sequence of $(\Uq,\Uq)$-modules
\beq
\label{shes}
{
 %(\hat{R}_q W \otimes W -  W \otimes W \hat{R}_q) =
0 \rightarrow (S_{2}(V^\ast) \otimes S^{1^2}(V)
\oplus S_{1^2}(V^\ast) \otimes S^{2}(V) ) 
%\oplus 
 %S_{(1^2)}(V^\ast) \otimes S^{(2)}(V)) 
\rightarrow \C(q)\la W \ra \rightarrow \C[GL_q(V)] \rightarrow 0} 
\eeq
whose polarization  yields the short exact sequence
of $(\hr,\hr)$-modules
\[
0 \rightarrow (S_{2} \otimes S^{1^2}  \oplus S_{1^2} \otimes S^{2}) 
\rightarrow \hr \otimes \hr \stackrel{p}{\rightarrow} \hr 
\cong \hr\otimes_{\hr} \hr 
\rightarrow 0 \ .
\]
The projection $p: \hr \otimes \hr {\rightarrow} \hr $
acts by $ p( T^\alpha \otimes T_{\beta})=
 T^\alpha \otimes_{\hr} T_\beta $. 
Clearly 
$p(S_{2} \otimes S^{1^2})=S_{2} \otimes_{\hh} S^{1^2}=0
%= S_{1^2} \otimes_{\hh} S^{2}
=p(S_{1^2} \otimes S^{2})$.
The one-dimensional $\hh$-bimodules $S_{2} \otimes S^{1^2}$
and $S_{1^2} \otimes S^{2}$ can be visualized by the braid
diagrams \footnote{The span  %of both bimodules 
$S_{1^2} \otimes S^{2} \oplus S_{1^2} \otimes S^{2}$
is equivalent to the Hecke relations eqs (\ref{pol},\ref{heckemove}), see also \cite{oleg}.}
$$
\vcenter{
\xy 0;
/r1pc/:
 ,{\vcross\vcross-}
\endxy}
\; = \;
\vcenter{
\xy 0;
/r1pc/:
 (0,0)*{}; (0,-2)*{}; **\dir{-};
(1,0)*{}; (1,-2)*{}; **\dir{-};
\endxy}
\;+
q \,\,\,\, \vcenter{
\xy 0;
/r1pc/:
%(0,0)*{}; (0,-1)*{}; **\dir{-};
,{\vcross-};
\endxy
\xy 0;
/r1pc/:
(0,0)*{}; (0,-1)*{}; **\dir{-};
%(1,0)*{}; (1,-1)*{}; **\dir{-};
(1,0)*{}; (1,-1)*{}; **\dir{-};
\endxy} \,\,
 - q^{-1} \,\,\;
\vcenter{
\xy 0;
/r1pc/:
(0,0)*{}; (0,-1)*{}; **\dir{-};
,{\vcross-};
%(0,-1.5)*{}; (0,-2)*{}; **\dir{-};
(1,0)*{}; (1,-1)*{}; **\dir{-};
%(1,-1.5)*{}; (1,-2)*{}; **\dir{-};
\endxy}
\qquad \mbox{and} \qquad
\vcenter{
\xy 0;
/r1pc/:
 ,{\vcross\vcross-}
\endxy}
\; = \;
\vcenter{
\xy 0;
/r1pc/:
 (0,0)*{}; (0,-2)*{}; **\dir{-};
(1,0)*{}; (1,-2)*{}; **\dir{-};
\endxy}
\;-
q^{-1} \,\,\,\,\vcenter{
\xy 0;
/r1pc/:
%(0,0)*{}; (0,-1)*{}; **\dir{-};
,{\vcross-};
\endxy
\xy 0;
/r1pc/:
(0,0)*{}; (0,-1)*{}; **\dir{-};
%(1,0)*{}; (1,-1)*{}; **\dir{-};
(1,0)*{}; (1,-1)*{}; **\dir{-};
\endxy}
 +q \,\,\, \;
\vcenter{
\xy 0;
/r1pc/:
(0,0)*{}; (0,-1)*{}; **\dir{-};
,{\vcross-};
%(0,-1.5)*{}; (0,-2)*{}; **\dir{-};
(1,0)*{}; (1,-1)*{}; **\dir{-};
%(1,-1.5)*{}; (1,-2)*{}; **\dir{-};
\endxy} \, \,  .
%\vcenter{
%\xy 0;
%/r1pc/:
%(0,0)*{}; (0,-1)*{}; **\dir{-};
%,{\vcross-};
%(1,0)*{}; (1,-1)*{}; **\dir{-};
%\endxy} =
%\vcenter{
%\xy 0;
%/r1pc/:
%,{\vcross-};
%\endxy
%\xy 0;
%/r1pc/:
%(0,0)*{}; (0,-1)*{}; **\dir{-};
%(1,0)*{}; (1,-1)*{}; **\dir{-};
%\endxy}
 %\quad .
$$

%Suppose we consider $\hr$
%as a right reducible $\hr$-module. By
%fixing the upper index
%$\alpha=12 \ldots r$
%we identify the row elements $(\hr)^{12 \ldots r}_\beta$ with 
%the Hecke generators $T_\beta$. On the other side,
 %by fixing
%$\beta=12 \ldots r$ we get the 
%column elements $(\hr)_{12 \ldots r}^\alpha$ identified with 
%$T^\alpha$ yielding a basis of the left $\hr$-module.
%One has $T^\alpha=T_{\alpha^{-1}}$
%and $$T^\alpha_\beta=T^{\alpha}\otimes_{\hr} T_\beta=
%1 \!\!1 \otimes_{\hr} T_{\alpha^{-1}}T_\beta $$
%which can be understood as choosing a representative in
%the coset
%\beq
%\hr \cong \hr \otimes_{\hr}\hr
%\eeq
%The idea of this note is to consider
\section{Quantum Diagonal algebra $\dia$}

\begin{definition} The quantum diagonal algebra $\dia$  is the subalgebra  of  the  quantum matrix algebra $\C[GL_q(V)]$
generated by the elements $x^1_1, x^2_2, \ldots, x^d_d$.
\end{definition}
The restriction of the commutative ring  $\C[GL(V)]$ to the subring of the diagonal matrix elements $x^i_i$ yields $\C[\T]$,  the {\it commutative} functions on the torus $\T$.
We now derive the relations in the restriction
$\dia$ to the diagonal of the  noncommutative ring
$\C[GL_q(V)]$.

\begin{lemma} Let  $\T^\ast= \sum_{i=1}^d \C x^i_i= W^{\Delta}$
be the  span of the diagonal generators in $\C[GL_q(V)]$.
%The FRT-module $\hec(W)$  generates the ideal of the $\C[GL_q(V)]$-relations
%\beq
%\label{abcc}
%\xymatrix{0 \ar@{->}[r] & (\hec(W) ) \ar@{^{(}->}[r]
%&T(W)
 %\ar@{>>}[r]& \C[GL_q(V)] \ar@{->}[r]  & 0} \ .
%\eeq
%The restriction of
%$\C[GL_q(V)]$ to  $\T^{\ast \otimes 3}$
%yields the short exact sequence
%$\T^{\ast \otimes 3} \subset W^{\otimes 3}$
%\[
%\xymatrix{0 \ar@{->}[r] & \pK(\T) \ar@{^{(}->}[r]
%&\T^{\ast\otimes 3}
 %\ar@{>>}[r]& \C[GL_q(V)]^\Delta_3 \ar@{->}[r]  & 0} \ .
%\]
The subspace ${\pK(\T)}\subset \T^{\ast \otimes 3}$ of the cubic  relations
of the quantum diagonal algebra $\dia$ is  generated by
the pseudo-plactic relations  (\ref{ppK}).
Its dimension is 
$\dim {\pK(\T)} =
\left( \ba{c} d \\ 3\ea  \right) +
2 \left( \ba{c} d \\ 2\ea \right)$ 
where  
$d=\dim V$.

%\beq
%\label{ppK}
%\ba{rclcrcl}
         %{\pK(\T)}^{i_1i_2}_{i_3}  &:=&
   %[  [ x_{i_1}^{i_1},x_{i_3}^{i_3} ], x_{i_2}^{i_2}]  &\qquad&\mbox{with}&\qquad
%&i_1<i_2<i_3
   %\\ [4pt]
 %{\pK(\T)}^{i_1i_1}_{i_2}  &:=&
        %[[ x_{i_1}^{i_1}, x_{i_2}^{i_2} ], x_{i_1}^{i_1} ]_{q^2}
%& & \mbox{with}& & i_1<i_2
 %\\[4pt]
 %{\pK(\T)}^{i_1i_2}_{i_2}  &:=&
  %[  x_{i_2}^{i_2}, [ x_{i_1}^{i_1}, x_{i_2}^{i_2}
%]]_{{q^2}}               && \mbox{with} & & i_1<i_2
         %\ea \ . \eeq  	
\end{lemma}

{\bf Proof. }
By expanding %$x^{ijk}_{ijk}$
$x^i_i x^j_j x^k_k=:x^{ijk}_{ijk}$ in  a basis $x_{ijk} :=x^1_i x^2_j x^3_k$ %=:\left[{}^{123}_{ijk}\right]$
 of
monomials in
$\C[GL_q(V)]_3$
one gets the linear dependences between the elements
in $\dia$, the nontrivial
solutions of the equation $ x^\sigma_\sigma c^\sigma=0 \ .$

Expanding the diagonal monomials
  $x^{ijk}_{ijk}$ in the basis $x_{ijk}$ yields the matrix
	%(up to a transposition)	
	\beq
\label{systd}
\left(
\ba{c}
x^{123}_{123} \\
x^{132}_{132} \\
x^{213}_{213}\\
x^{231}_{231}\\
x^{312}_{312} \\
x^{321}_{321}
\ea
\right)
%\left(
%\ba{c}
%\left[{}^{123}_{123}\right] \\[4pt]
%\left[{}^{132}_{132}\right] \\[4pt]
%\left[{}^{213}_{213}\right]\\[4pt]
%\left[{}^{231}_{231}\right]\\[4pt]
%\left[{}^{312}_{312}\right] \\[4pt]
%\left[{}^{321}_{321}\right]
%\ea
%\right)
=
\left(\ba{cccccc}
1 &0 & 0& 0& 0 & 0\\
1& \omega &0 &0&0&0
\\
1&0 & \omega & 0 &0&0
\\
1&0 &\omega & 0&0 & \omega
\\
1& \omega &0 & 0&0 & \omega \\
1& \omega & \omega & \omega^2 & \omega^2 & \omega^3+\omega
\ea
\right)
\left(
\ba{c}
x_{123} \\
x_{132} \\
x_{213}\\
x_{231}\\
x_{312} \\
x_{321}
\ea
\right) \qquad \qquad \omega=q-q^{-1} \ .
\eeq
Thus one has a  transformation  with a singular matrix,
%leading to $c_\sigma x^\sigma_\sigma=0$
$
x^\sigma_\sigma = x_\rho M_{\sigma}^{\rho}$ with $\rho, \sigma\in S_3
%\qquad
$.
Since the monomials $x_\rho$ are a true basis in
$1^3$-graded part of $\C[GL_q(V)]_3$ we have
\[
 x^\sigma_\sigma c^\sigma=0 \quad \Leftrightarrow \quad
M_{\sigma}^{\rho}  c^\sigma=0 \qquad\quad  {\rm Ker}{M} = \pK^{12}_3
=\C \, [  [ x_{1}^{1},x_{3}^{3} ], x_{2}^{2}] \ ,
\]
where the combination of %the sencond, third, forth and fifth
columns (rows in eq (\ref{systd})) of
$M_{\sigma}^{\rho}$ yields the nontrivial  solution
$
c^{132}=c^{231}=-c^{312}=-c^{213}=1
$ and $c^{123}=0=c^{321}$ .
%Hence
%the 1-dimensional kernel of $M_{\sigma}^{\rho}$ is 
%$\lambda[  [ x_{1}^{1},x_{3}^{3} ], x_{2}^{2}]$, $\lambda \in \C$.
%The kernel of the matrix $M$ spans the ideal of the relations
%${\pK(\T)}^{i_1i_2}_{i_3} $.

For  $x^\sigma_\sigma\in \dia$ of weight $(2,1)$  we have a basis $x_{112}:=
x^{112}_{112}$ and $x_{121}:=x_{121}^{112}$. Hence
by expanding the diagonal elements 
we get a $3\times 2$ matrix $M$
\[
\left(\ba{c}
x_{112}^{112}  \\
x_{121}^{121}\\
x_{211}^{211}
\ea
\right) 
=
\left(
\begin{array}{cc}
 1 & 0 \\
 1 & \omega \\
 1 & (q^2+1) \omega   \\
\end{array}
\right)
\left(\ba{c}
x_{112} \\
x_{121}
\ea
\right)  \qquad \quad {\rm Ker} M = \pK^{11}_2
=\C \, [[ x_{1}^{1}, x_{2}^{2} ], x_{1}^{1} ]_{q^2}
\]
%The 1-dimensional kernel yields
%$[[ x_{1}^{1}, x_{2}^{2} ], x_{1}^{1} ]_{q^2}$. 
Similarly from $x^\sigma_\sigma\in \dia$ of weight $(1,2)$ one obtains
${\rm Ker} M = \pK^{12}_2$. $\Box$
% =\C \,[  x_{2}^{2}, [ x_{1}^{1}, x_{2}^{2}  ]]_{{q^2}}$. %We are done.
% $\Box$

{\bf Remark.}
The origin of the kernel of $M$
 is the existence of  two ``homotopic'' expressions of
 maximal element
$x_{321}$ % of weight $(1,1,1)$ 
whose difference %yields 
$[  [ x_{1}^{1},x_{3}^{3} ], x_{2}^{2}]\omega^{-1} \in{\pK}^{1,2}_{3}$ %in eq. (\ref{ppK})
\beq
\label{ambig}
x_{321}^{123}=
x^1_3x^3_1 x_{2}^{2}= [ x_{3}^{3},x_{1}^{1} ] x_{2}^{2} \omega^{-1}
\qquad
x_{321}^{123}= x_{2}^{2} x^1_3 x^3_1= x_{2}^{2} [ x_{3}^{3},x_{1}^{1} ]  \omega^{-1}
 \ .
\eeq
In the same vein, the maximal element 
of weight $(2,1)$, that is, $x_{211}^{112}$
while expressed in  $\dia$
can be written in two different ways
\beq
x^{112}_{211}=q^{-1}x^{121}_{211}=
q^{-1} [x^2_2,x^1_1]x^1_1 \omega^{-1}
\qquad 
x^{112}_{211}=qx^{112}_{121}=
 q x^1_1 [x^2_2,x^1_1] \omega^{-1}
\eeq
and the difference of two expressions
yields  $-\omega^{-1} q^{-1}[ [x^1_1,x^2_2],x^1_1]_{q^2} \in {\pK}^{1,1}_{2}$. 
% in eq. (\ref{ppK})
Similarly
 the maximal element 
of weight $(1,2)$, {\it i.e.}, $x_{221}^{122}$ forks and
leads to ${\pK}^{1,2}_{2}$.

%\beq
%x^{122}_{221}=qx^{122}_{212}=q [x^2_2,x^1_1]x^2_2 \omega^{-1}
%\qquad 
%x^{122}_{221}=q^{-1}x^{212}_{221}=
 %q^{-1} x^2_2 [x^2_2,x^1_1] \omega^{-1}
%\eeq

\section{A functorial way to 
the diagonal algebra $\C[GL_q(V)]^\Delta$}
\label{pre-pl}
\begin{definition} Let  
$\C[GL_q(V)]^\Delta$ be the diagonal algebra generated by 
$x^i_i\in \C[GL_q(V)]$.
The diagonal Hecke algebra
 ${\h}^\Delta= \bigoplus_r \hr^\Delta$ is the polarization  
of $\C[GL_q(V)]^\Delta$ %= \bigoplus_r \C[GL_q(V)]_r^\Delta$
\[
\C[GL_q(V)]^\Delta_r =\sum_{A\in \{1, \ldots,d \}^r} \C(q) x^A_A \xymatrix{\ar@{~>}[r]&}  
\hr^\Delta = \sum_{\alpha\in \S_r} \C(q) T^{\alpha}_{\alpha} \ .
\]
\end{definition}

We fix a partition of the unit in $1\!\!1_{\hhh} $ by orthogonal idempotents
\[
1\!\!1_{\hhh} = e_3 + e_{21}^+ + e_{21}^- + e_{1^3}
\]
%these are 
%the $q$-symmetrizer $e_3$, the $q$-antisymmetrizer 
%$e_{1^3}$ and
 the two idempotents  $e_{21}^\pm$
corresponding to the two Standard Young Tableaux with shape $\lambda=21$,
$\dim S_{21}=2$.
 The idempotent ${e}^+_{21}$ is a deformation of
 the Eulerian idempotent $e_3^{[1]}$ considered by Jean-Louis Loday \cite{eulerien}.
 Eulerian idempotents $e_n^{[1]}$ split
the Harisson homology from the Hochschild homology.
\begin{lemma}{(\cite{LP2})}
The splitting of the central idempotent
$E_{21}=e^+_{21} + e^-_{21} $
%splits  
into two minimal idempotents
$e^+_{21}e^-_{21}=0$  is uniquely chosen
 by the  $\theta$-multiplication eigenvalues
\beq
\label{idempY}
\theta e^{\pm}_{21}= \pm e^{\pm}_{21} =  e^{\pm}_{21} \theta
% \tau \qquad \qquad \tau
\qquad
\qquad \theta= T_{s_1}T_{s_2}T_{s_1} \ .
\eeq
%where $\theta= T_{s_1}T_{s_2}T_{s_1}$
 %stands for the maximal element in $\hhh$.
These minimal idempotents  are
polynomials  of $T_{\sigma} \in \hhh$ (for details see \cite{LP2})
\beqa
  \label{idemp}
%e(q):=
e^{\pm}_{21}&=&
\frac{1}{[3]}\left(T_{123} -\frac{1}{2} (T_{231}\pm T_{213}\pm T_{132}+T_{312})
\pm T_{ 321}\right) \\
\nonumber
&+&\frac{\omega}{2[3]} \left( T_{213}\mp T_{312}  \mp T_{231} + T_{132}\right) \qquad \omega = q-q^{-1} \ .
\eeqa
The projector
${e}^-_{21}$ is obtained
from   $e^+_{21}$ by the involution ${T}_\sigma \rightarrow (-1)^\sigma T_\sigma$, $q \rightarrow q^{-1}$.
\end{lemma}
%{\bf Proof.} See the appendix.

 %The two irreducible isomorphic right (left) $\hhh$-modules are denoted by
%$
 %S^\pm_{21} =e^{\pm}_{21} \hhh $ $
%(S_\pm^{21} = \hhh e^{\pm}_{21}) 
%$
We prove in the appendix the following  important lemma
\begin{lemma} Let us denote by  $\mathfrak{L}_q^\pm(W)$
the $\Uq$-bimodule 
$$\mathfrak{L}_q^\pm(W)=
e_{21}^\pm W^{\otimes 3}e_{21}^\mp= 
S^\pm_{21}(V^\ast) \otimes S_\mp^{21}(V) \ .$$
The relations of the quantum
 pseudo-plactic algebra $\PP_q(\T)$
are the image of  the restriction of $\mathfrak{L}_q^\pm(W)$  to the diagonal $\T^\ast$
%of a $\Uq$-bimodule 
%$\mathfrak{L}_q^+(W)=e_{21}^+ W^{\otimes 3}e_{21}^-$
\beq
\pK(\T)=\mathfrak{L}_q^+(W)|_{\T^\ast} =\mathfrak{L}_q^-(W)|_{\T^\ast}  \ .
\label{funct}
\eeq
%$$ \pK(\T) = \pK(W)|_{\T^\ast}  $$
\label{brute}
\end{lemma}
%{\bf Proof.} See the appendix.
%By abuse of notation we will write $e^\mp_{21}$
%for the image in the $\hhh$-representation $\pi(e^\mp_{21})$,
%{\it i.e.}, a polynomial of the matrices $\pi(T_{s_1})=(\hat{R}_{q})_{12}$ and $\pi(T_{s_2})=(\hat{R}_{q})_{23}$
%in ${\rm{End}}(\Vt 3)$. Then the matrices $\pi(T_{s_i})$
%commute with
%quantum matrix entries %in $W^{\otimes 3}$
 %hence
%the orthogonality  $e^+_{21}e^-_{21}=0$ implies
 %$p(\mathfrak{L}_q^\pm(W))=p(e^+_{21}W^{\otimes 3} e^-_{21})=0$.

\begin{definition}  The pre-plactic algebra $\PP_q$
is the graded algebra 
$\PP_q = \bigoplus_{r\geq 0} \PP_q(r)$ with degrees
given by the quotient
\[
\PP_q (r) \cong (\hr \otimes \hr)^\Delta/ (\pK)_r
\]
where  $(\pK)_r$ stays for  the degree $r$ of the ideal $(\pK)$  generated
by the polarization  of the  pseudo-plactic 
relations ${\pK}(\T)$ cf. eq. (\ref{funct})
\[
\pK :=
[[13]2]:=
\widetilde{T}^{132}_{\hspace{2pt} 132}-\widetilde{T}^{312}_{\hspace{2pt} 312}
-\widetilde{T}^{213}_{\hspace{2pt}213}+\widetilde{T}^{231}_{\hspace{2pt} 231}\qquad \qquad \widetilde{T}^\alpha_\alpha :=T^\alpha\otimes T_\alpha \ .
\]
\end{definition}
It is clear that the  pre-plactic algebra 
$\PP_q \subset \h =\bigoplus_{r\geq 0}\hr$
is  the polarization of the quantum pseudo-plactic algebra
$\PP_q(\T)$. We obtain now the key result;

\begin{theorem} 
The diagonal Hecke algebra  
$(\h)^\Delta =\bigoplus_{r\geq 0}(\hr)^\Delta
$ is isomorphic
to the pre-plactic algebra
 $\PP_q= \bigoplus_{r\geq 0}\PP_q(r)$
\[(\h)^\Delta \cong \PP_q \ .\]
\label{pp}
\end{theorem}

%By construction $p(\pK)=0$. 

{\bf Proof.}
The $r!$ elements
$\widetilde{T}^\alpha_\alpha :=T^\alpha\otimes T_\alpha $
freely generate the % involution invariant 
diagonal space
$$(\hr \otimes \hr)^\Delta := \bigoplus_{\alpha\in \S_r}
\C(q) \widetilde{T}^\alpha_\alpha \ .$$ 
We are going to show now
that  $\pK$
is the unique combination in $(\hhh \otimes \hhh)^\Delta$
 which is projected out by 
$p: \hhh \otimes \hhh \rightarrow \hhh \otimes_{\hhh} \hhh$, 
by a pictorial proof
$$p(\pK)={T}^{132}_{\hspace{2pt} 132}-{T}^{312}_{\hspace{2pt} 312}
-{T}^{213}_{\hspace{2pt}213}+{T}^{231}_{\hspace{2pt} 231}=0 
$$

We attach to each
generator in $(\calh_3(q)\otimes\calh_3(q))^\Delta$
its braid using coset notation (\ref{coset1})
putting the left factor in the upper
half-plane (above the horizon) and
right factor under the horizon
\[
\ba{cccccccc}
T^{132}_{\hspace{2pt}132} &=& T^{132}_{\hspace{2pt}123} \otimes
T^{123}_{\hspace{2pt}132}& =& T_{s_2}\otimes T_{s_2} &=&
\vcenter{
\xy 0;
/r1pc/:
{\vcross};
(0,0)*{}; (0,-1)*{}; **\dir{-};
\endxy
\xy 0;
/r1pc/:
{\vcross};
(0,0)*{}; (0,-1)*{}; **\dir{-};
\endxy
}
\\[14pt]
T^{312}_{\hspace{2pt}312} &=& T^{312}_{\hspace{2pt}123} \otimes
T^{123}_{\hspace{2pt}312} &=& T_{s_2 s_1}\otimes T_{s_1 s_2}
&=&
\vcenter{
\xy 0;
/r1pc/:
{\vcross};
(0,0)*{}; (0,-1)*{}; **\dir{-};
\endxy
\xy 0;
/r1pc/:
{\vcross};
(3,0)*{}; (3,-1)*{}; **\dir{-};
\endxy
\xy 0;
/r1pc/:
{\vcross};
(3,0)*{}; (3,-1)*{}; **\dir{-};
\endxy
\xy 0;
/r1pc/:
{\vcross};
(0,0)*{}; (0,-1)*{}; **\dir{-};
\endxy
}
 \\[14pt]
T^{231}_{\hspace{2pt}231} &=& T^{231}_{\hspace{2pt}123} \otimes
T^{123}_{\hspace{2pt}231}&= &T_{s_1s_2}\otimes T_{s_2s_1}&=&\quad
 \vcenter{
\xy 0;
/r1pc/:
{\vcross};
(3,0)*{}; (3,-1)*{}; **\dir{-};
\endxy
\xy 0;
/r1pc/:
{\vcross};
(0,0)*{}; (0,-1)*{}; **\dir{-};
\endxy
\xy 0;
/r1pc/:
{\vcross};
(0,0)*{}; (0,-1)*{}; **\dir{-};
\endxy
\xy 0;
/r1pc/:
{\vcross};
(3,0)*{}; (3,-1)*{}; **\dir{-};
\endxy
}
	\\[14pt]
	T^{213}_{\hspace{2pt}213} &=& T^{213}_{\hspace{2pt}123} \otimes
T^{123}_{\hspace{2pt}213} &=& T_{s_1}\otimes T_{s_1}&=&
\vcenter{
\xy 0;
/r1pc/:
{\vcross};
(3,0)*{}; (3,-1)*{}; **\dir{-};
\endxy
\xy 0;
/r1pc/:
{\vcross};
(3,0)*{}; (3,-1)*{}; **\dir{-};
\endxy
} 	
\ea
\]
These braids are symmetric with respect to the horizon,
the left factor being the braid indexed by the inverse permutation of the right.
The projection $p$ is gluing the upper and the lower
braids allowing generators to flow across 
 the tensor product
$$p(\widetilde{T}^{\alpha}_{\alpha}) =
T_{\alpha^{-1}}\otimes_{\hhh} T_{\alpha}=
1\!\!1 \otimes_{\hhh} T_{\alpha^{-1}} T_{\alpha}$$
and allows to reduce the number of crossings, by reducing the word written with braid generators.
For instance, in the Hecke relation 
$T_{s_1}^2=1\!\!1+ (q-q^{-1})T_{s_1}$, 
the ``bubble" $(T_{s_1})^2$ being reduced  to braids in $\h$ 
with
lower length number of crossings
\beq
\label{heckemove}
\vcenter{
\xy 0;
/r1pc/:
 ,{\vcross\vcross-}
\endxy}
\; = \;
\vcenter{
\xy 0;
/r1pc/:
 (0,0)*{}; (0,-2)*{}; **\dir{-};
(1,0)*{}; (1,-2)*{}; **\dir{-};
\endxy}
\;+ \omega\;
\vcenter{
\xy 0;
/r1pc/:
(0,0)*{}; (0,-.5)*{}; **\dir{-};
,{\vcross-};
(0,-1.5)*{}; (0,-2)*{}; **\dir{-};
(1,0)*{}; (1,-.5)*{}; **\dir{-};
(1,-1.5)*{}; (1,-2)*{}; **\dir{-};
\endxy} 
%%%%%%%%%%%%%%%%%%%%%%%%%%%%%%
%%%%%%%%%%%%%%%%%%%%%%%%%%
\qquad  \qquad
\vcenter{
\xy 0;
/r1pc/:
(0,0)*{}; (0,-.5)*{}; **\dir{-};
,{\vcross-};
(0,-1.5)*{}; (0,-2)*{}; **\dir{-};
(1,0)*{}; (1,-.5)*{}; **\dir{-};
(1,-1.5)*{}; (1,-2)*{}; **\dir{-};
\endxy} :=
\vcenter{
\xy 0;
/r1pc/:
(0,0)*{}; (0,-1)*{}; **\dir{-};
,{\vcross-};
(1,0)*{}; (1,-1)*{}; **\dir{-};
\endxy} =
\vcenter{
\xy 0;
/r1pc/:
,{\vcross-};
\endxy
\xy 0;
/r1pc/:
(0,0)*{}; (0,-1)*{}; **\dir{-};
(1,0)*{}; (1,-1)*{}; **\dir{-};
\endxy}
\qquad  \omega=  (q- q^{-1})
\eeq
a move that we are referring as Hecke move.

The ``standardized" pseudo-plactic relation is pictorially represented by
sum of ``diagonal'' diagrams
\[
T^{132}_{\hspace{2pt} 132}-T^{312}_{\hspace{2pt} 312}
-T^{213}_{\hspace{2pt}213}+T^{231}_{\hspace{2pt} 231}=
\vcenter{
\xy 0;
/r1pc/:
{\vcross};
(0,0)*{}; (0,-1)*{}; **\dir{-};
\endxy
\xy 0;
/r1pc/:
{\vcross};
(0,0)*{}; (0,-1)*{}; **\dir{-};
\endxy
} \quad  - \; \quad
\vcenter{
\xy 0;
/r1pc/:
{\vcross};
(0,0)*{}; (0,-1)*{}; **\dir{-};
\endxy
\xy 0;
/r1pc/:
{\vcross};
(3,0)*{}; (3,-1)*{}; **\dir{-};
\endxy
\xy 0;
/r1pc/:
{\vcross};
(3,0)*{}; (3,-1)*{}; **\dir{-};
\endxy
\xy 0;
/r1pc/:
{\vcross};
(0,0)*{}; (0,-1)*{}; **\dir{-};
\endxy
}
\quad  + \; \quad
\vcenter{
\xy 0;
/r1pc/:
{\vcross};
(3,0)*{}; (3,-1)*{}; **\dir{-};
\endxy
\xy 0;
/r1pc/:
{\vcross};
(0,0)*{}; (0,-1)*{}; **\dir{-};
\endxy
\xy 0;
/r1pc/:
{\vcross};
(0,0)*{}; (0,-1)*{}; **\dir{-};
\endxy
\xy 0;
/r1pc/:
{\vcross};
(3,0)*{}; (3,-1)*{}; **\dir{-};
\endxy
}
\quad -\quad \vcenter{
\xy 0;
/r1pc/:
{\vcross};
(3,0)*{}; (3,-1)*{}; **\dir{-};
\endxy
\xy 0;
/r1pc/:
{\vcross};
(3,0)*{}; (3,-1)*{}; **\dir{-};
\endxy
} 	=
\]
which we simplify by reducing 
 terms $T^{312}_{\hspace{2pt} 312}$ and $T^{231}_{\hspace{2pt} 231}$
\[
\vcenter{
\xy 0;
/r1pc/:
{\vcross};
(0,0)*{}; (0,-1)*{}; **\dir{-};
\endxy
\xy 0;
/r1pc/:
{\vcross};
(0,0)*{}; (0,-1)*{}; **\dir{-};
\endxy
} \quad  - \; \quad
%%%%%%%%%%%%%%%%%%%%%%%%%%%%%%%%%%%%%%%%%
\vcenter{
\xy 0;
/r1pc/:
{\vcross};
(0,0)*{}; (0,-1)*{}; **\dir{-};
\endxy
\xy 0;
/r1pc/:
(1,0)*{}; (1,-1)*{}; **\dir{-};
(2,0)*{}; (2,-1)*{}; **\dir{-};
(3,0)*{}; (3,-1)*{}; **\dir{-};
\endxy
\xy 0;
/r1pc/:
(1,0)*{}; (1,-1)*{}; **\dir{-};
(2,0)*{}; (2,-1)*{}; **\dir{-};
(3,0)*{}; (3,-1)*{}; **\dir{-};
\endxy
\xy 0;
/r1pc/:
{\vcross};
(0,0)*{}; (0,-1)*{}; **\dir{-};
\endxy
}
%%%%%%%%%%%%%%%%%%%%%%%%%%%%%%
\quad  - \; \omega \,
\vcenter{
\xy 0;
/r1pc/:
{\vcross};
(0,0)*{}; (0,-1)*{}; **\dir{-};
\endxy
\xy 0;
/r1pc/:
(0,0)*{}; (0,-.5)*{}; **\dir{-};
,{\vcross-};
(0,-1.5)*{}; (0,-2)*{}; **\dir{-};
(1,0)*{}; (1,-.5)*{}; **\dir{-};
(1,-1.5)*{}; (1,-2)*{}; **\dir{-};
(2,0)*{}; (2,-2)*{}; **\dir{-};
\endxy
\xy 0;
/r1pc/:
{\vcross};
(0,0)*{}; (0,-1)*{}; **\dir{-};
\endxy
}
%%%%%%%%%%%%%%%%%%%%%%%%%%%%%%%%%%%%%
\quad  + \; \omega \,
\vcenter{
\xy 0;
/r1pc/:
{\vcross};
(3,0)*{}; (3,-1)*{}; **\dir{-};
\endxy
\xy 0;
/r1pc/:
(0,0)*{}; (0,-2)*{}; **\dir{-};
(1,0)*{}; (1,-.5)*{}; **\dir{-};
,{\vcross-};
(1,-1.5)*{}; (1,-2)*{}; **\dir{-};
(2,0)*{}; (2,-.5)*{}; **\dir{-};
(2,-1.5)*{}; (2,-2)*{}; **\dir{-};
\endxy
\xy 0;
/r1pc/:
{\vcross};
(3,0)*{}; (3,-1)*{}; **\dir{-};
\endxy
}
%%%%%%%%%%%%%%%%%%%%%5555
\quad  + \; \quad
\vcenter{
\xy 0;
/r1pc/:
{\vcross};
(3,0)*{}; (3,-1)*{}; **\dir{-};
\endxy
\xy 0;
/r1pc/:
(0,0)*{}; (0,-1)*{}; **\dir{-};
(1,0)*{}; (1,-1)*{}; **\dir{-};
(2,0)*{}; (2,-1)*{}; **\dir{-};
\endxy
\xy 0;
/r1pc/:
(0,0)*{}; (0,-1)*{}; **\dir{-};
(1,0)*{}; (1,-1)*{}; **\dir{-};
(2,0)*{}; (2,-1)*{}; **\dir{-};
\endxy
\xy 0;
/r1pc/:
{\vcross};
(3,0)*{}; (3,-1)*{}; **\dir{-};
\endxy
}
%%%%%%%%%%%%%%%%%%%%%%%%%%%%%%%
\quad -\quad \vcenter{
\xy 0;
/r1pc/:
{\vcross};
(3,0)*{}; (3,-1)*{}; **\dir{-};
\endxy
\xy 0;
/r1pc/:
{\vcross};
(3,0)*{}; (3,-1)*{}; **\dir{-};
\endxy
} \ .
\]
The rightmost and leftmost terms with bubbles 
are homotopic and cancel hence
$$
T^{132}_{\hspace{2pt} 132}-T^{312}_{\hspace{2pt} 312}
-T^{213}_{\hspace{2pt}213}+T^{231}_{\hspace{2pt} 231}=
-(q-q^{-1}) \,\left( \quad \vcenter{
\xy 0;
/r1pc/:
{\vcross};
(0,0)*{}; (0,-1)*{}; **\dir{-};
\endxy
\xy 0;
/r1pc/:
,{\vcross};
(3,0)*{}; (3,-1)*{}; **\dir{-};
\endxy
\xy 0;
/r1pc/:
{\vcross};
(0,0)*{}; (0,-1)*{}; **\dir{-};
\endxy }
\quad
 - \quad \;
\vcenter{
\xy 0;
/r1pc/:
,{\vcross};
(3,0)*{}; (3,-1)*{}; **\dir{-};
\endxy
\xy 0;
/r1pc/:
{\vcross};
(0,0)*{}; (0,-1)*{}; **\dir{-};
\endxy
\xy 0;
/r1pc/:
,{\vcross};
(3,0)*{}; (3,-1)*{}; **\dir{-};
\endxy }
\quad \right) = 0
$$
and we got the transparent result: the
pre-plactic relation $\pK$ %(standardized pseudo-Knuth relation)
is equivalent to the braid relation.
The reduced word $321$ has two representatives and
the two  braids cancel 
$$p(\pK)=\omega (T_{s_1}T_{s_2}T_{s_1}- T_{s_2}T_{s_1}T_{s_2})=0 \ .$$
\begin{lemma} All symmetric elements in $(\h \otimes \h)^\Delta$ with vanishing projection in $\hdia$ belong to
the ideal generated by the braid relation
$(\pK)$,
\[(\pK) = \{x\in (\h \otimes \h)^\Delta|  p(x)=0 \} \ . \]
\label{exact}
\end{lemma}
\vspace{-8pt}
{\bf Proof of the lemma.}
Reducing a generator $x$ in $\hdia \subset \h$ to its minimal length
is a  ''normal ordering'' such that
 we can't apply the Hecke moves 
the braid diagram of $x$ any more. Any non-zero element
in $\hdia$
is represented by a combination of reduced words,
but there is a remaining ''gauge freedom'',
$T_{s_is_{i+1}s_i}\sim T_{s_{i+1}s_{i}s_{i+1}} $ which is
% The remaining freedom for a word of minimal length is
the ``mutation" of the reduced word with respect to the braid relation.

Assume that in degree $n>3$ we have a relation $\hdia$  not generated by $\pK$.
 We conclude that there exists
 an independent relation between the reduced words in $\h$
 which is of degree higher than 3 which is a contradiction
since any reduced word can be brought to any other by a sequence
of braiding mutations $T_{s_is_{i+1}s_i}\sim T_{s_{i+1}s_{i}s_{i+1}} $.
The lemma is proven.

When restricted to the diagonal $\T$
the sequence  of $(\Uq,\Uq)$-modules cf. eq. (\ref{shes})
yields the sequence of spaces % $\Uq$-modules
\beq
\label{she}
{
 0 \rightarrow (\pK(\T))_r 
\rightarrow \C(q)\la \T \ra \rightarrow \C[GL_q(V)]^\Delta \rightarrow 0}  \ .
\eeq
These spaces are stable under the quantum Weyl action \cite{LS}
which lives in a  completion 
$\widehat{U_q{\mathfrak{gl}}}(V)$.

The diagonal restriction
$(\h{\otimes } \h)^{\Delta}:= \bigoplus_{r\geq 0}\left(\sum_{\alpha\in \S_r} \C \widetilde{T}^\alpha_{\, \, \alpha}\right)$
  of the $\h$-bimodule $\h\otimes \h$ 
% $(\h{\otimes } \h)^{sym}:= \bigoplus_{r\geq 0}\left(\bigoplus_{\alpha\in \S_r} \C T^\alpha_{\, \, \alpha}\right)$ 
is a $\h$-module with a left action
$$T_\rho \cdot \widetilde{T}_{\,\,\sigma}^\sigma:=
T_{\rho^{-1}} T_{\sigma^{-1}}\otimes T_\sigma T_{{\rho}}
%, \qquad\sigma, \rho \in \S_r \ . 
\qquad (\mbox{similarly for} 
\quad
 {T}_{\,\,\sigma}^\sigma \in \hdia) \ .
$$
The polarization of this  sequence of $\widehat{U_q{\mathfrak{gl}}}(V)$-modules provides the sequence of $\hr$-modules
\beq
\label{ses}
{
 0 \rightarrow (\pK)_r  
\rightarrow (\hr \otimes \hr)^\Delta \stackrel{p}{\rightarrow}
 \hr^\Delta \rightarrow 0}  \ .
\eeq
 %The projection $p$ maps   
%$p(\widetilde{T}^\alpha_\alpha) =
%T^\alpha\otimes_{\hr} T_\alpha$ hence
%$p (\pK)=0$.
According to lemma
\ref{exact} the sequence of $\hr$-modules cf. (\ref{ses}) is exact for all $r\geq 0$. 
The exactness  implies the isomorphism $(\h)^\Delta \cong \PP_q$. 
The theorem is proven. $\Box$

Theorem \ref{pp} implies the conjecture of Daniel Krob and Jean-Yves Thibon \cite{KT}. % is now follows. 
%a corollary of .

\begin{corollary}
The diagonal algebra $\dia$ is isomorphic to
the quantum pseudo-plactic algebra
\[
\dia \cong \PP_q(\T) \ . %:= \C (q)\la \T \ra / (\pK(\T)) 
\]
\end{corollary}
The isomorphism $\PP_q(\T) \cong \dia$ holds true if and only if the  sequence of $\widehat{U_q{\mathfrak{gl}}}(V)$-modules (\ref{she}) is exact
\[
%\PP_q(\T) =
\dia \cong \C (q)\la \T \ra / (\pK(\T)) \ .
\]
By functoriality the exactness of the  sequence 
of $\widehat{U_q{\mathfrak{gl}}}(V)$-modules 
(\ref{she})  follows from the exactness of
the sequence of $\hr$-modules (\ref{ses}).
The latter exactness is due 
the isomorphism $\PP_q\cong \hdia$ (Theorem \ref{pp}).

\begin{acknowledgement} It is my pleasure to
thank Peter Dalakov, Tekin Dereli, Vladimir Dobrev, Michel Dubois-Violette, G\'erard Duchamp, 
Ludmil Hadjiivanov, Nikolay Nikolov, Petko Nikolov, Oleg Ogievetsky and
Ivan Todorov   for their encouraging interest in that work and many enlightening discussions.
This work has been supported  by the Bulgarian National Science Fund research grant DN 18/3 and in part by the
TUBITAK 2221 program.

%If you want to include acknowledgments of assistance and the like at
%the end of an individual chapter please use the
%\verb|acknowledgement| environment -- it will automatically render
%Springer's preferred layout.
\end{acknowledgement}
\section*{Appendix}
\addcontentsline{toc}{section}{Appendix}

{\bf Proof of Lemma \ref{brute}.}
By abuse of notation we will write $e^\mp_{21}$
for the image in the $\hhh$-representation $\pi(e^\mp_{21})$,
{\it i.e.}, a polynomial of the matrices $\pi(T_{s_1})=(\hat{R}_{q})_{12}$ and $\pi(T_{s_2})=(\hat{R}_{q})_{23}$
in ${\rm{End}}(\Vt 3)$. The matrices $\pi(T_{s_i})$ %, $i=1,2$
%$\hat{R}_{12}$ and $\hat{R}_{23}$ 
commute with
quantum matrix elements
%\[
%(\hat{R}_{q})_{i i+1} W\otimes W\otimes W= W\otimes W\otimes W (\hat{R}_{q})_{i i+1}
%\qquad i=1,2
%\]
\[
(\hat{R}_{q})_{12} W^{\otimes 3}= W^{\otimes 3} (\hat{R}_{q})_{12} \qquad 
(\hat{R}_{q})_{23} W^{\otimes 3}= W^{\otimes 3} (\hat{R}_{q})_{23}
%\qquad i=1,2
\]
thus the orthogonality  $e^+_{21}e^-_{21}=0$ implies
 $p(\mathfrak{L}_q^\pm(W))=p(e^+_{21}W^{\otimes 3} e^-_{21})=0$.
%thus belongs to the ideal $(\pK(\T))_3$.

The proof is by brute force, a direct check
using the $n^3 \times n^3$ matrix $[e^\pm_{21}]^{i_1 i_2 i_3}_{j_1 j_2 j_3} $. 

The grading of the Drinfeld-Jimbo matrix eq. (\ref{DJ}) implies
that the matrices in ${\rm{End}}(\Vt 3)$ have zero entries
 $[e^\pm_{21}]^{i_1 i_2 i_3}_{j_1 j_2 j_3} = 0$ if
%\qquad \mbox{if} \qquad 
$\{i_1 i_2 i_3\}\neq \{j_1 j_2 j_3\}$
% \quad  \mbox{as multisets}$$  
as  multisets
 thus it is enough to restrict our attention
to $\dim V=3$.
The matrix of the idempotent $[e^\pm_{21}]$ has  $6\times 6$ blocks for 3 different indices
and  $3\times 3$ blocks for 2 different indices.
The $6\times 6$ blocks are indexed by
$\sigma, \rho \in\{{123}, {132}, {213}, {231}, {312}, {321}\}$
$$
[e^\pm_{21}]^{\sigma}_{\rho}=
\frac{1}{[3]}\left(
\begin{array}{cccccc}
 1 & \frac{\omega \mp 1}{2} &
   \frac{\omega \mp 1}{2} &
     -\frac{1\pm \omega }{2} &
    -\frac{1\pm \omega }{2}   &\pm 1 \\
 \frac{\omega \mp 1}{2} &
	 \frac{\omega^2 \mp \omega +2}{2}
	&
    -\frac{1 \pm \omega}{2}   &\pm 1 &
   \mp \frac{\omega^2+1}{2} &
   -\frac{1 \mp \omega}{2}\\
 \frac{\omega \mp 1}{2} &
     -\frac{1 \pm \omega}{2}   &
  \frac{\omega^2 \mp \omega +2}{2} &
   \mp \frac{\omega^2+1}{2} &
   \pm 1 &  -\frac{1 \mp \omega}{2}
   \\
 -\frac{1\pm \omega }{2}& \pm 1 &
   \mp \frac{\omega^2+1}{2}&
    \frac{\omega^2 \pm \omega +2}{2} &
   -\frac{1\mp \omega}{2} &
    -\frac{\omega\pm 1}{2}  \\
  -\frac{1\pm \omega }{2} &
    \mp \frac{\omega^2+1}{2} &
  \pm 1 & -\frac{1\mp \omega}{2}&
    \frac{\omega^2 \pm \omega +2}{2} &
    -\frac{\omega\pm 1}{2}  \\
 \pm 1 & -\frac{1 \mp \omega}{2}&
  -\frac{1 \mp \omega}{2} &
    -\frac{\omega\pm 1}{2}  &
    -\frac{\omega \pm 1}{2}  & 1 \\
\end{array}
\right)
$$
%$\sigma, \rho \in\{{i_1i_2i_3}, {i_1i_3i_2}, {i_2i_1i_3}, {i_2i_3i_1}, {i_3i_1i_2}, {i_3i_2i_1}\}$
while the $3\times 3$ blocks are indexed by
$\lambda, \mu \in \{112,121,211\}$ or $ \{122,212,221\}$.
\[
[e^\pm_{21}]^\lambda_\mu=\frac{1}{2(q +q^{-1}\pm 1)}
\left(
\begin{array}{ccc}
 q & -1\mp q & \pm 1 \\
 -1 \mp q &  q\pm 2+ q^{-1} & -1 \mp q^{-1} \\
 \pm 1 &- 1 \mp q^{-1} & q^{-1} \\
\end{array}
\right)
%\qquad
%\left(
%\begin{array}{ccc}
 %q & q-1 & -1 \\
 %q-1 & q-2+\frac{1}{q} & \frac{1}{q}-1 \\
 %-1 & \frac{1}{q}-1 & \frac{1}{q} \\
%\end{array}
%\right)
\]

Applying the Einstein summation convention
 over repeating indices
$k$ and $l$ but not on $i$'s
we get the sum
\[
%\sum_{k_1,k_2,k_3,l_1,l_2,l_3}
{\mathfrak{L}_q^\pm}(W)^{i_1 i_2 i_3}_{i_1 i_2 i_3} =
% [e^\pm(q) W\otimes W\otimes W e^\mp(q)]^{i_1i_2i_3}_{i_1i_2i_3}=
[e^\pm(q)]^{i_1i_2i_3}_{k_1k_2k _3}
x^{k_1}_{l_1}  x^{k_2}_{l_2} x^{k_3}_{l_3}
[e^\mp(q)]^{l_1l_2l_3}_{i_1i_2i_3} \qquad \mbox{no summation on $i_1$, $i_2$, $i_3$} \ .
\]
For  multi-indices $1=i_1< i_2 = i_3=2$ of weight $(1,2)$
we get
%It is enough to look at $1=i_1<i_2=2$
\[
{\mathfrak{L}_q^\pm}(W)^{ 1 2 2}_{ 1 2 2} =
[e^\pm(q)]^{122}_{abc}
x^{a}_{i}  x^{b}_{j} x^{c}_{k}
[e^\mp(q)]^{ijk}_{122} =\frac{[2] }{4 \omega[3]}[x^2_2,[x^1_1,x^2_2]]_{q^2}
\in {\pK(\T)}^{12}_{2}
\]
where the off-diagonal terms have  replaced by the substitutions
\beq
\ba{ccccc}
x^{122}_{221}\rightarrow q x^{122}_{212} &&
x^{221}_{122}\rightarrow q x^{212}_{122} &&
x^{122}_{212}\rightarrow [x^2_2,x^1_1]x^2_2 /\omega
\\x^{221}_{212}\rightarrow x^{212}_{221}&&
x^{212}_{122}\rightarrow x^{122}_{212}
&&
x^{212}_{221}\rightarrow x^2_2[x^2_2,x^1_1]/\omega
\ea \ .
\eeq
%leading to
%\[
%{\mathfrak{L}_q^\pm}(W)^{2 1 2}_{2 1 2}
%=- \frac{[2] }{4 q[3]}[x^2_2,[x^1_1,x^2_2]]_{q^2}
%\sim {\pK(\T)}^{12}_{2} \ .
%\]
By similar substitutions   for indices of weight 
$(2,1)$ we get
 ${\mathfrak{L}_q^\pm}(W)^{ 1  12}_{ 1  12}
\in {\pK(\T)}^{11}_{2} \ .$
%\[
%{\mathfrak{L}_q^\pm}(W)^{1 1 2}_{1 1 2} =
%[e^\pm(q)]^{112}_{abc}
%x^{a}_{i}  x^{b}_{j} x^{c}_{k}
%[e^\mp(q)]^{ijk}_{112} =\frac{[2] }{4 q[3]}[x^2_2,[x^1_1,x^2_2]]_{q^2}
%\sim {\pK(\T)}^{11}_{2}\ .
%\]
%Since the quantum group  relations
%depend only on the relative order we have obtained
%\[
%\ba{ccccccc}
%[\mathfrak{L}^\pm_q(\T)]^{i_1i_2i_1}_{i_1i_2i_1}
%&=&
 %\frac{[2] }{4 [3]\omega}
  %[[ x_{i_1}^{i_1}, x_{i_2}^{i_2} ], x_{i_1}^{i_1} ]_{q^2}
%&\sim&	
%{\pK(\T)}^{i_1i_1}_{i_2}  &&i_1<i_2\\[4pt]
%[\mathfrak{L}^\pm_q(\T)]^{i_2i_1i_2}_{i_2i_1i_2}
%&=&
 %\frac{[2] }{4 [3]\omega}
  %[x_{i_2}^{i_2},[ x_{i_1}^{i_1}, x_{i_2}^{i_2} ] ]_{q^2}
%&\sim&	
%{\pK(\T)}^{i_1i_2}_{i_2}  &&i_1<i_2
%\ea \ .
%\]
%Hence we recover all pseudo-Knuth relatations
%with repeating indices in (\ref{ppK}) through projectors $e^{\pm}(q)$.

For  multi-indices %$i_1< i_2 <i_3$ 
of weight $(1,1,1)$
%we get We now show how to
we ''diagonalize"
the expression %of weight $(1,1,1)$
  $${\mathfrak{L}_q^\pm}(W)^{1 2 3}_{1 2 3}=
[e^\pm(q)]^{123}_{\sigma}
x^{\sigma}_{\rho}
[e^\mp(q)]^{\rho}_{123}\qquad \sigma,\rho\in S_3$$
by rewriting it in diagonal monomials
$x^\sigma_\sigma:=x^i_i x^j_j x^k_k\in (\T^\ast)^{\otimes 3}$
where $\sigma\in S_3 $ is the word $ijk$ of the permutation
$\sigma(1)=i, \sigma(2)=j, \sigma(3)=k$. %  which belong to $\pK(\T)$.
 The restriction to the diagonal subalgebra $\pK(\T)$
can be done by
the substitutions
\[
\ba{ccccc}
x^{213}_{321} \rightarrow x^{123}_{ 231} +
  \omega x^{213}_{ 231} &&
	x^{132}_{ 321} \rightarrow
 x^{123}_{312} + \omega x^{132}_{312} &&
x^{231}_{312}\rightarrow
 x^{123}_{231} + \omega x^{213}_{231} \\
x^{213}_{312}\rightarrow x^{123}_{132}   +\omega x^{123}_{312}
&&
x^{312}_{231}\rightarrow x^{123}_{312}
  + \omega x^{123}_{321} &&
	x^{321}_{213}\rightarrow x^{123}_{312}
	+\omega x^{231}_{213}\\
	x^{312}_{213} \rightarrow x^{132}_{123}+\omega x^{132}_{213}
	&&
	x^{231}_{132}\rightarrow x^{213}_{123} +\omega x^{123}_{312}&&
	x^{132}_{231}\rightarrow x^{123}_{213}+ \omega x^{123}_{231}\\
	x^{132}_{213} \rightarrow x^{123}_{231} &&
	x^{321}_{132}%\rightarrow x^{213}_{312}
	\rightarrow x^{123}_{132}   +\omega x^{123}_{312} &&
	x^{213}_{132}\rightarrow x^{123}_{312}\\
	x^{321}_{123}\rightarrow x^{231}_{213}&&
	x^{312}_{123}\rightarrow x^{123}_{231} &&
	x^{231}_{123}\rightarrow x^{123}_{312}
	%x[1, 2] ** x[3, 1] ** x[2, 3] -> x[1, 2] ** x[2, 3] ** x[3, 1]
	%x[3, 1] ** x[2, 3] ** x[1, 2] -> x[2, 3] ** x[1, 2] ** x[3, 1],
%x[2, 1] ** x[1, 3] ** x[3, 2] -> x[1, 3] ** x[2, 1] ** x[3, 2],
%x[3, 1] ** x[2, 2] ** x[1, 3] -> x[2, 2] ** x[3, 1] ** x[1, 3],
%x[3, 1] ** x[1, 2] ** x[2, 3] -> x[1, 2] ** x[2, 3] ** x[3, 1],
%x[2, 1] ** x[3, 2] ** x[1, 3] -> x[1, 3] ** x[2, 1] ** x[3, 2]
\\
\ea
\]
where the terms in  RHS are either lowest in lexicographical order
, {\it i.e., } $x^{123}_{abc}$
or with one  rightmost (leftmost)  diagonal entry , $x^i_k x^k_i x^j_j$ and
$x^j_j x^i_k x^k_i$. In the  latter case one can
systematically replace   the off-diagonal terms by commutators
of diagonal terms
\[
x^i_k x^k_i x^j_j \rightarrow {[x^i_i,x^k_k]x^j_j}/ {\omega}
\qquad \qquad
x^j_j x^i_k x^k_i  \rightarrow {x^j_j[x^i_i,x^k_k]}/ {\omega}
\ .
\]
For instance, the terms $x_{213}:=x^{123}_{213}$ and $x_{132}:=x^{123}_{132}$ which are also lowest in lexicographical order are rewritten as
\beq
\label{offdia}
x^{123}_{213}= [x^2_2,x^1_1]x^3_3 / {\omega}
\qquad
x^{123}_{132}= x^1_1[x^3_3,x^2_2] / {\omega} \ .
\eeq
After imposing all substitution above we are left with only
three non-diagonal terms $x^{123}_{231}$, $x^{123}_{312}$ and
$x^{123}_{321}$.

The last equation in the system (\ref{systd}) is the only one containing
 $x^{123}_{231}$ and $x^{123}_{312}$
\[
x^{321}_{312}-
x^{123}_{123} - \omega x^{123}_{213} -\omega
x^{123}_{132} -(\omega^3 + \omega)x^{123}_{321}
= \omega^2(x^{123}_{231}+x^{123}_{312})
\] hence it provides an obstruction
of an expression to be reducible to a sum of
$x^i_i x^j_j x^k_k\in (\T^\ast)^{\otimes 3}$:
it is clear that a sum
$\sum_{\sigma\in S_3}d_\sigma x^{123}_\sigma \in (\T^\ast)^{\otimes 3}$
if and only if $d_{231}=d_{312}$.

The direct check shows that indeed in the sum
${\mathfrak{L}_q^\pm}(W)^{1 2 3}_{1 2 3}$
the coefficient of  $x^{123}_{231}$ is the same as
the coefficient of $x^{123}_{312}$
thus the substitution\footnote{We have used
the substitution (\ref{offdia}) to eliminate off-diagonal
terms $- \omega x^{123}_{213} -\omega
x^{123}_{132}$}
\[
x^{123}_{312} \rightarrow -x^{123}_{231} +
\frac{1}{\omega^2}\left\{ x^{321}_{312}-
x^{123}_{123} - [x^2_2,x^1_1]x^3_3- x^1_1[x^3_3,x^2_2]
 -(\omega^3 + \omega)x^{123}_{321} \right\}
\]
will cancel term $x^{123}_{312}$ with $x^{123}_{231}$.
Finally we can  eliminate the last off-diagonal term
$x^{123}_{321}$ by one of the two possible ways
(\ref{ambig}) or a combination thereof.
We are going to choose the ``gauge''
\[
x^{123}_{321} \rightarrow \omega^{-1} [x^3_3,x^1_1]x^2_2 \ .
\]
By direct calculation % (that we checked by computer algebra)
 after replacing all
off-diagonal terms by the above substitution
 we recover the pseudo-Knuth relations $\pK(\T)$ for three different indices in eq. (\ref{ppK})
\[
[\mathfrak{L}^\pm_q(\T)]^{i_1i_2i_3}_{i_1i_2i_3}
=
-\frac{\omega^3 \mp \omega^2 \mp 2 }{2\omega[3]^2} %\frac{[3]+2}{ 2[3]^2}
   [  [ x_{i_1}^{i_1},x_{i_3}^{i_3} ], x_{i_2}^{i_2}]
\in	{\pK(\T)}^{i_1i_2}_{i_3} \qquad i_1<i_2<i_3
 \ .
\]

%$$\hec(W) \cong S_{(1,1)}(V^{\ast})\otimes_{\calh_2(q)} S^{(2)}(V)
%\cong S_{(2)}(V^{\ast})\otimes_{\calh_2(q)} S^{(1,1)}(V) \ .$$

%\input{referenc}

\begin{thebibliography}{99.}%
% Use the following syntax and markup for your references if
% the subject of your book is from the field
% "Mathematics, Physics, Statistics, Computer Science"
%
% Contribution

%\bigskip

% Use the following (APS) syntax and markup for your references if
% the subject of your book is from the field
% "Mathematics, Physics, Statistics, Computer Science"
%
% Online Document

\bibitem{D-VP} M. Dubois-Violette and T. Popov,
         Homogeneous algebras, statistics and combinatorics.
          {\it Lett. Math. Phys. \/ \bf 61} (2002), 159-170.

\bibitem{FRT} L. D.
Faddeev, N. Yu. Reshetikhin, and L. A. Takhtajan. "Quantization of Lie groups and Lie algebras." Algebraic analysis. Academic Press, 1988. 129-139.
\bibitem{KT}
D. Krob and J.-Y. Thibon, Noncommutative symmetric functions IV: Quantum linear groups and Hecke algebras at $q= 0$.{ \it Journal of Algebraic Combinatorics \/ \bf 6}(1997),  339-376.
\bibitem{Plaxique} A. Lascoux, B. Leclerc, and J.-Y. Thibon, The plactic monoid. {\em Algebraic Combinatoric on Words} (2002), 10pp.

\bibitem{LS}
Levendorski, S. Z., and Ya S. Soibelman. "Some applications of the quantum Weyl groups." Journal of Geometry and Physics 7.2 (1990): 241-254.

\bibitem{eulerien} J.-L. Loday,
S\'erie de Hausdorff, idempotents eul\'eriens et alg\`ebres de Hopf.
{\em Exposition. Math. 12 }(1994), 165-178.

\bibitem{LP2} {J.-L. Loday and T. Popov},
 Parastatistics Algebra, Young Tableaux and Super Plactic Monoid.
 {\em	International Journal of Geometric Methods in Modern Physics\/ \bf 5} (2008),  1295-1314.

\bibitem{LP3} J.-L. Loday, T. Popov.
Hopf Structures on Standard Young Tableaux.
{\em Proceedings % of the International Workshop
 "Lie Theory and Its Applicatioins"},%' in Physics", Varna  2009},
  Vl. Dobrev ed.,
  AIP conference series {Vol. \bf 1243}(2010), 265--275.
	
	\bibitem{MR} C. Malvenuto and C. Reutenauer, Duality
between Quasi-Symmetric Functions and the Solomon Descent Algebra.
{ \em Journal of Algebra {\bf 177}} (1995), 967--982.

 \bibitem{Manin} Yu. I. Manin, Quantum groups and noncommutative geometry.
Univer-sit\'e de Montr\'eal, Centre de Recherches Math\'ematiques, Montr\'eal, QC,1988.

\bibitem{oleg}
O.
Ogievetsky, Uses of quantum spaces. {\em Contemporary Math. {\bf 294} }
(2002), 161-232.


	\bibitem{PR} S.  Poirier and C. Reutenauer,  Alg\`ebres
de Hopf de tableaux. {\em Ann. Sci. Math. Qu\' ebec \/ \bf 19 } (1995),
79-90.

\bibitem{TP}
T. Popov,  "Quantum Plactic and Pseudo-Plactic Algebras." International Workshop on Lie Theory and Its Applications in Physics. Springer, Singapore, 2015.

\bibitem{TP2}
T. Popov,  "Pre-Plactic Algebra and Snakes." arXiv preprint arXiv:1711.06253 (2017).

%\bibitem{phys-online} J. Dod, in \textit{The Dictionary of Substances and Their Effects}, Royal Society of Chemistry. (Available via DIALOG, 1999),
%\url{http://www.rsc.org/dose/title of subordinate document. Cited 15
%Jan 1999}

\end{thebibliography}
\end{document}